\newcommand{\Real}{\mathbb{R}}
\newcommand{\Gel}{\mathbb{G}}
\newcommand{\Kel}{\mathbb{K}}

\newcommand{\Ratio}{\mathbb{Z}}
\newcommand{\Natural}{\mathbb{N}}

\newcommand{\caL}{\mathcal L}
\newcommand{\caG}{\mathcal G}
\newcommand{\caH}{\mathcal H}
\newcommand{\caD}{\mathcal D}
\newcommand{\caS}{\mathcal S}
\newcommand{\caM}{\mathcal M}
\newcommand{\caP}{\mathcal P}

\newcommand{\Diss}{\frak{D}}

\newcommand{\mb}{\mathbf{m}}
\newcommand{\Kb}{\mathbf{K}}

\def\dive{\mathop{\rm div}}
\def\diag{\mathop{\rm diag}}
\def\span{\mathop{\rm span}}
\def\meas{\mathop{\rm meas}}
\newcommand{\dif}{\mathrm{d}}
\newcommand{\ljump}{[\![}
\newcommand{\rjump}{]\!]}

\documentclass[preprint, 12pt]{elsarticle}
\usepackage{latexsym}
\usepackage{mhchem}
\usepackage{amsfonts,amssymb,amsmath, amsbsy}
\usepackage{graphicx}

\usepackage{tikz}
\usepackage{pgf}
\usepackage{pgffor}
\usepgfmodule{shapes}
\usepgfmodule{plot}
\usetikzlibrary{decorations}
\usetikzlibrary{arrows}
\usetikzlibrary{snakes}

\usepackage{color}
\usepackage{tikz}
\usepackage{calc}

\newlength{\yellownotewidth}
\newlength{\yellownoteheight}
\journal{Journal of Total Rejection}

\begin{document}
\begin{frontmatter}
\title{A brief overview of existence results and decay time estimates for a mathematical modeling of scintillating crystals}
\author{Fabrizio Dav\'{\i}\footnote{e-mail: davi@univpm.it}}
\address{DICEA \& ICRYS\\ Universit\'a Politecnica delle Marche, via Brecce Bianche, 60131 Ancona, Italy\\
on leave at IMT - Lucca, Italy}

\begin{abstract} 
Inorganic scintillating crystals can be modelled as continua with microstructure. For rigid and isothermal crystals the evolution of charge carriers becomes in this way described by a reaction-diffusion-drift equation coupled with the Poisson equation of electrostatic. Here we give a survey of the available existence and asymptotic decays results for the resulting boundary value problem, the latter being a direct estimate of the scintillation decay time.  We also show how to recover various approximated models which encompass also the two most used phenomenological models for scintillators, namely the Kinetic and Diffusive ones. Also for these cases we show, whenever it is possible, which existence and asymptotic decays estimate results are known to date. 
\end{abstract}
\end{frontmatter}

\section{Introduction}

A scintillator crystal is a material which converts ionizing radiations into photons in the frequency range of visible light, hence its name. It acts as a true "wavelength shifter" and in such a role is used as radiation sensor into high-energy physics, in medical imaging and in security applications \cite{Lecoq}. The physics of scintillation, which is a  complex multi-scale phenomenon (\emph{see e.g.\/} \cite{VG14}) can be described within a continuum approach at three scales: at a \emph{Microscopic} scale the incoming energy $E$ generates a population of charged energy carriers which moves in straight directions for few nanometer \cite{JA07} and whose density $N=N(E)$ can be found by the means of approximated solutions of the Bethe-Bloch equation \cite{SE52}-\cite{ULMA10}.\footnote{ In \cite{D1} we show how the density of excitation carriers $N$ induced by an ionizing energy $E$ which hits the crystal at a given point $x^{*}$, can be obtained by the means of a suitable rescaling to the mesoscopic scale of the approximate solutions of the Bethe-Bloch equation along an elementary cylindrical track: in such a way $N=N(E)$ maintains informations on both the initial energy $E$ and the material properties of the crystal.} These energy carriers wander and migrate within a greater region either generating other energy carriers or recombining with emission of photons $h\nu$. In the process some energy is lost and a scintillator is a material in which such a loss reduces the frequency of the incoming ionizing energy to that of visible light. We call this scale the \emph{Mesoscopic} scale: for $\caP$ the region occupied by the crystal we denote $\Omega\subset\caP$ the \emph{mesoscopic volume} in which the recombination of charge carriers into photons takes place (Fig. 1). Finally the light rays propagate within the crystal at a \emph{Macroscopic} scale according to the laws of classical optics. 
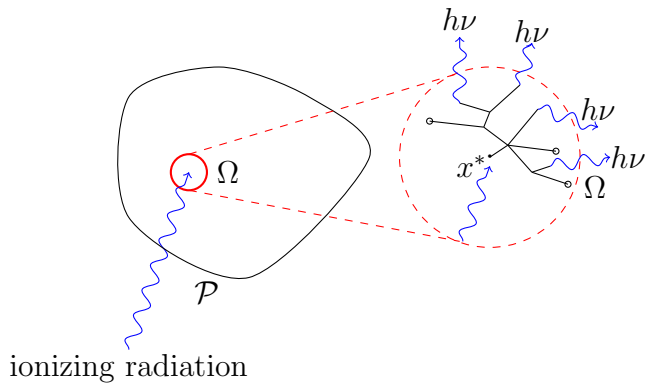
\begin{figure}
\begin{center}
\begin{tikzpicture}[scale=0.8]
\draw[color=red, thick] (6,3.75) circle (0.3cm);
\draw[black](6.3,1.75) node {${\caP}$};
\draw[black](12.75,3.5) node {${\Omega}$};   
\draw[black](6.65,3.75) node {${\Omega}$};     
\draw[color=blue, ->, decorate, decoration=snake] (5,0.8)--(6,3.75);
\draw[black](5,0.5) node {ionizing radiation};  
  \draw[color=blue, ->, decorate, decoration=snake] (10.55,2.6)--(11,3.85);  
\draw plot[smooth cycle] coordinates{(5,3) (7,2) (9,4) (8,5) (6.5,5.5) (5,5)};         
\draw[color=red, dashed] (6,4.05)--(10.5,5.4);
\draw[color=red, dashed] (11,4) circle (1.5cm);
\draw[black] (11,4) node {$\cdot$};
\draw[color=red, dashed] (6,3.45)--(10.5,2.6);
\draw[black] (10.7, 3.85) node {$x^{*}$};
\draw[black](11,4)--(11.3,4.2);

\draw[black](11.3,4.2)--(10.9,4.5);
\draw[black](10.9,4.5)--(10,4.6);
\draw[color=black] (10,4.6) circle (0.05cm);

\draw[black](10.9,4.5)--(11,4.75);
\draw[black](11,4.75)--(10.5, 4.9);
\draw[color=blue, ->, decorate, decoration=snake] (10.5,4.9)--(10.5,6);
\draw[black] (10.5,6.3) node {$h\nu$};

\draw[black](11,4.75)--(11.5, 5.2);
\draw[color=blue, ->, decorate, decoration=snake] (11.5, 5.2)--(11.7,5.9);
\draw[black] (11.7,6.1) node {$h\nu$};

\draw[black](11.3,4.2)--(11.8,4.8);
\draw[color=blue, ->, decorate, decoration=snake] (11.8,4.8)--(12.8,4.5);
\draw[black] (12.8,4.8) node {$h\nu$};

\draw[black](11.3,4.2)--(11.7,3.75);

\draw[black](11.7,3.75)--(12,3.85);
\draw[color=blue, ->, decorate, decoration=snake] (12,3.85)--(13,4);
\draw[black] (13.3,4) node {$h\nu$};

\draw[black](11.7,3.75)--(12.3,3.55);
\draw[color=black] (12.3,3.55) circle (0.05cm);

\draw[black](11.3,4.2)--(12.1,4.1);
\draw[color=black] (12.1,4.1) circle (0.05cm);

\end{tikzpicture}
\end{center}
 \caption{The mesoscopic scintillation region $\Omega$.}
 \label{Davi_fig:01}
\end{figure}

Scintillation is a fast and dissipative phenomenon and there are two major physical parameters which are to be improved in  a scintillating material: (i) the \emph{decay time} $\tau_{d}$, which is the time required for scintillation emission to decrease to $e^{-1}$ of its maximum and is a measure of the scintillator resolution; (ii) the \emph{light yield} $LY$, which is the ratio between the collected light energy and the energy of the incoming ionizing radiation and which is a measure of the scintillator efficiency.

Scintillators can be modelled as Continua with microstructure \cite{D1}-\cite{D5} to arrive at a Reaction-Diffusion-Drift (RDD) equation for the energy carriers descriptors, coupled with the Poisson equation of electrostatic and both with Neumann-type boundary conditions. Further in \cite{D2}, by following the results obtained in \cite{CHJU17} and \cite{FK18}, we showed that for these equations it is possible to proof the global existence of renormalized and weak solutions and also how the decay time can be estimated explicitly in terms of these equations  constitutive parameters.

Here, by using  the results obtained into \cite{CHJU18}, we expand the results of \cite{D2} to show the existence of weak-strong renormalized solutions:  we further show also (whenever it is possible) existence and asymptotic decay results for the phenomenological models widely used in the literature, which can be obtained from our model by introducing suitable approximations as we did into \cite{D3}. 

We remark that scintillation is strongly affected by temperature and the crystal are linearly elastic deformable bodies, but here we limit our analysis to the isothermal case (the temperature being at most a parameter in the constitutive quantities) and to rigid crystals as in \cite{D3}. The effects of temperature are dealt with into \cite{D1} and \cite{D4} whereas an insight into the deformable case is provided in \cite{D5}. A complete treatment of these electromagnetical, thermal and mechanical interactions will be provided in a forthcoming book \cite{DAbook}.

The paper is organized as follows: in \S.\textbf{2} we give an overview of the model proposed into \cite{D1}-\cite{D5}, which leads to a RDD system: then we study the properties of the stationary solutions, which are important when we deal with the solutions asymptotic properties.  

In \S.\textbf{3} We extend with the help of \cite{CHJU18} results presented into \cite{D2} and we give results on existence of weak-strong  renormalized solutions and on asymptotic decay,  this last result being a direct estimate of the decay time. This result allows, for the first time, to estimate the decay time in terms of the parameters describing the physical properties of a given scintillating crystal and which appears explicitly in the equations we obtained.

Finally in \S.\textbf{4} first of all we put the coupled RDD and Poisson equations in an adimensional form by introducing characteristic length and time related to the scale of scintillation, as we did into \cite{D3}. The resulting boundary value problem depends on a set of three adimensional parameter related to diffusion, drift and recombination respectively. Then by choosing a suitable measure of smallness we obtain, to within higher-order terms in such a measure, different approximated models which represents different physically meaningful regimes for scintillators. In such a way we not only recover the most used phenomenological model in use, namely the \emph{Kinetic} and \emph{Diffusive} ones, but also show under which hypotheses scintillation could be described either by a Reaction-Diffusion or a Diffusion-Drift equation. Also for these approximated boundary value problems we give a survey of the existence and asymptotic decay results, either by taking them straight from the available literature or by adapting them to the specific context of scintillation as described by our model.

\section{The continuum model for isothermal and rigid scintillators: the boundary value problem}

\subsection{A reaction-diffusion-drift equation for scintillators}

The charge carriers with density $N$ generated at the point $x^{*}$ and time $t^{*}$ and whose dimension is (length)$^{-3}$, represent a population of carriers which can differ by their sign (\emph{e.g.} negative electrons, positive holes, neutral excitons and so on) and by their recombination mechanism. Accordingly we may differentiate these carriers by introducing an ordered array $n$, the \emph{charge carrier densities } vector:
\begin{equation}\label{Nvector}
n\equiv(n_{1}\,,n_{2}\,,\ldots n_{k})\,,\quad N=\sum_{j=1}^{k}n_{j}\,;
\end{equation}
clearly, the greater is $k$, the finer will be our description of the phenomena: the simplest non-trivial choice is $k=2$ as in \cite{LGW11}, where $n_{1}$ represents the electrons population (which is equal to the holes population) and $n_{2}$ represents the population of excitons, which are bounded electron-hole pairs evolving together. In the various phenomenological models for scintillation thus far proposed we may have $m=3$ as in \cite{WGLU11} and \cite{BM12}, $m=7$ as in \cite{GBD15} whereas in \cite{VA08} we have $m\geq 11$. The charge carrier densities are related to $N$ by:
\begin{equation}
n_{j}=\alpha_{j}N\,,\quad\alpha_{j}\geq0\,,\quad j=1\,,2\,,\ldots\,,k\,,\quad\sum_{j=1}^{k}\alpha_{j}=1\,;
\end{equation}
the set $\{\alpha_{1}\,,\alpha_{2}\,,\ldots\,,\alpha_{k}\}$ depends on the specific scintillator, on the initial energy $E$ and on the initial data (\emph{cf.\/} the discussion in \cite{BM09}).

We identify the charge carriers densities at $(x^{*}\,,t^{*})$ with fields defined on the whole $\Omega\times[0\,,\tau)$ and whose regularity will shall made precise when needed:
\begin{equation}
\Omega\times[0\,,\tau)\ni(x\,,t)\mapsto n_{j}(x\,,t)\geq 0\,,\quad j=1\,,2\,,\ldots k\,;
\end{equation}
without loss of generality we may assume that the mesoscopic volume $\Omega$ is a ball of a still unprescribed radius $R$, centered at $x^{*}$.

We  assume accordingly  the $k-$dimensional field, the \emph{excitation carrier densities} vector
\begin{equation}
\Omega\times[0\,,\tau)\ni(x\,,t)\mapsto n(x\,,t)\in\caM\,,\quad\caM=\Real_{+}^{k}\cup\{0\}\,,
\end{equation}
as the main state variable in our description of scintillation. In view of (\ref{Nvector}), since $n_{j}\geq 0$, $j=1\,,2\,,\ldots\,,k$, we can view the extension of $N$ to a field on $\Omega$ as the  $L^{1}(\Omega)$ norm of $n$:
\begin{equation}\label{L1norm}
\|n\|_{L^{1}(\Omega)}=\int_{\Omega}\sum_{j=1}^{k}|n_{j}|=\int_{\Omega} N\,.
\end{equation}

Let $e$ be the elementary charge, then the excitation carrier vector induces a free charge density within the scintillation volume $\Omega$
\begin{equation}
\rho_{f}=ez\cdot n\,,\quad\mbox{ in }\Omega\,,
\end{equation}
with $z=(z_{1}\,,z_{2}\,,\ldots z_{k})$, $z_{j}\in\Ratio$, $j=1\,,2\,,\ldots k$ the \emph{charge vector}. By the Maxwell-Lorentz equations in absence of magnetic fields \cite{WA79} these free charges induce a \emph{local electric potential} $(x\,,t)\mapsto\varphi(x\,,t)$ which, for a given time $t\in[0\,,\tau)$, is the solution of the Poisson equation of electrostatic with associated Neumann boundary conditions:
\begin{eqnarray}\label{laplace}
&-\epsilon\Delta\varphi=\chi|_{{}_{\Omega}}ez\cdot n\,,\quad\mbox{ in }\Real^{3}\times[0\,,\tau)\,,\\
&\ljump\nabla\varphi\rjump\cdot\mb=0\,,\quad\mbox{ on }\partial\Omega\times[0\,,\tau)\,,\nonumber
\end{eqnarray}
where
\begin{equation}
\chi|_{{}_{\Omega}}=
\begin{cases}
1\,,\quad\mbox{ in }\Omega\,,\\
\\
0\,,\quad\mbox{ in }(\Real^{3}/\Omega)\,;
\end{cases}
\end{equation}
here $\epsilon$ is the permittivity of the crystal (which at this stage we assume isotropic or at most cubic),\footnote{For anisotropic crystals equations (\ref{laplace}) become:
\begin{eqnarray}\label{laplaceaniso}
&-\epsilon_{o}\dive\Kb[\nabla\varphi]=\rho^{*}\,,\quad\mbox{ in }\Real^{3}\times[0\,,\tau)\,,\\
&\ljump\Kb\nabla\varphi\rjump\cdot\mb=0\,,\quad\mbox{ on }\partial\Omega\times[0\,,\tau)\,,\nonumber
\end{eqnarray}
where $\epsilon_{o}$ is the vacuum permittivity and $\Kb$ is the symmetric and positive definite \emph{permittivity tensor}; relations (\ref{electroenergy})  changes accordingly with  $\epsilon_{o}\Kb[\nabla\varphi]\cdot\nabla\varphi$ in place of $\epsilon\|\nabla\varphi\|^{2}$} and $\mb$ is the outward unit normal to $\partial\Omega$. 

We notice that in (\ref{laplace}) we do not take into account either bound charges or external charges since we are mainly interested into $X-$ or $\gamma-$rays which have zero charge, whereas for $\alpha-$ and $\beta-$rays, which have respectively positive and negative charges, an external charge contribution $q^{*}$ should be added to (\ref{laplace}).

We further remark that an implicit way to select the radius of $\Omega$ s indeed that at its boundary the Neumann condition holds, \emph{i.e.} there is no electric field outflow trough $\partial\Omega$. 

The boundary-value problem (\ref{laplace}) admits an unique (up to a constant) weak solution $\varphi\in H^{1}(\Omega)$ provided the \emph{total charge} is  conserved:
\begin{equation}\label{conservationcharge}
Q(t)=e\int_{\Omega}z\cdot n=0\,,\quad\forall t\in[0\,,\tau)\,,
\end{equation}
annd the constant can be conveniently determined if we set $\bar{\varphi}=0, $\footnote{For any given integrable  function $f$ we shall denote with $\bar{f}$ its mean value:
\begin{equation*}
\bar{f}=\frac{1}{\meas(\Omega)}\int_{\Omega}f\,.
\end{equation*}} thus making the solution unique.

Following \cite{AGH02}, we assume that the \emph{electrostatic free-energy} associated to the electric potential has, besides the classical conservative term which depends on $n$ by the means of (\ref{laplace}), a dissipative term $F(n)$ of entropic nature which also depends on the excitation carrier vector:
\begin{equation}\label{electroenergy}
U(n)=\frac{1}{2}\int_{\Real^{3}}\epsilon\|\nabla\varphi(n)\|^{2}-\theta\int_{\Omega}F(n)\,,
\end{equation}
where $\theta>0$ is the (fixed) absolute temperature; in \cite{D4} we showed how to such an energy it can be associated an electrostatic self-power:
\begin{equation}
\dot{U}=\int_{\Omega}s\cdot\dot{n}\,,
\end{equation}
where the elements of the array $s\equiv(s_{1}\,,s_{2}\,,\ldots s_{k})$ represent the \emph{scintillation potentials} associated to the various charge carriers:
\begin{equation}
s_{j}(n)=ez_{j}\varphi(n)-\theta\frac{\partial F}{\partial n_{j}}(n)\,,\quad j=1\,,2\,,\ldots\,,k\,.
\end{equation}
We shall call $s$, with an abuse of terminology, the \emph{scintillation potentials vector}:
\begin{equation}\label{scintpotential}
s(n)=ez\varphi(n)-\theta\frac{\partial F}{\partial n}(n)\,,
\end{equation}
and it is easy to show that:
\begin{equation}\label{frechetU}
s=\caD U\,,
\end{equation}
where $\caD$ denotes the Frechet  derivative.

In \cite{D1} and the related papers \cite{D3}, \cite{D4} we showed how, by using a continuum with microstructure approach and the classical dissipation inequality, the evolution equation for the charge carriers in scintillators can be represented as:
\begin{eqnarray}\label{BVP1}
&\dive S(n)[\nabla s]-H(n)s=\dot{n}\,,\quad\mbox{ in }\Omega\times[0\,,\tau)\,,\\
&S(n)[\nabla s]\cdot\mb=0\,,\quad\mbox{ on }\partial\Omega\times[0\,,\tau)\,,\nonumber
\end{eqnarray}
where $S(n)$ and $H(n)$ are two definite-positive $k\times k$ matrices and
\begin{equation}\label{nablascint}
\nabla s=ez\otimes\nabla\varphi-\theta\frac{\partial^{2} F}{\partial n^{2}}\nabla n\,.
\end{equation}
The boundary-value problem (\ref{BVP1}) has an equivalent variational formulation which leads to a gradient-flow type problem (\emph{vid. e.g.\/} \cite{AGS05} and for more recent results \cite{MRS12}):
\begin{equation}\label{gradflow1}
\dot{n}=-\caD\Psi(n\,,s)\,,\quad\Psi(n\,,s)=\frac{1}{2}\int_{\Omega}S(n)[\nabla s]\cdot\nabla s+H(n)s\cdot s\,,
\end{equation}
where the \emph{conjugate dissipation} functional $\Psi(n\,,s)$ is related to the thermodynamical dissipation $\Diss$ by
\begin{equation}\label{dissipation1}
\Diss=2\Psi\,.
\end{equation}
The gradient-flow classical structure is 
\begin{equation}
\Gel(n)\dot{n}=-\caD U\,,
\end{equation}
and $U$ is the driving functional: starting with \cite{MIE11}, provided there exists $\Kel(n)=\Gel^{-1}(n)$, the following formulation was proposed and successfully used into \emph{e.g.\/} \cite{GMI13}-\cite{MAMI18}:
\begin{equation}\label{gradflow2}
\dot{n}=-\Kel(n)[\caD U]\,,
\end{equation}
where the Onsager structure $\Kel$ ca be splitted additively into different contributions. For instance, in our case we may set:
\begin{equation}
\Kel=\Kel_{D}+\Kel_{R}\,,\quad \Kel_{D}[\,\cdot\,]=-\dive S(n)\nabla[\,\cdot\,]\,,\quad \Kel_{R}[\,\cdot\,]=H[\,\cdot\,]\,.
\end{equation}

Furthermore, in \cite{D1} it is show that, provided we identify the entropic term with the Gibbs entropy, 
\begin{equation}\label{Gibbs}
F(n)=-k_{B}\sum_{j=1}^{k}n_{j}(\log(\frac{n_{j}}{c_{j}})-1)\,,
\end{equation}
with $c_{j}$, $j=1\,,2\,,\ldots\,,k$, normalizing constants and $k_{B}$ the Boltzmann constant then, by using (\ref{Gibbs}) into (\ref{scintpotential}) and (\ref{nablascint}), we arrive at a Reaction-Diffusion-Drift equation for the evolution of charge carriers in scintillators with Neumann boundary conditions and initial data:
\begin{eqnarray}\label{RDD}
&\dive(D[\nabla n]+M[N(n)z\otimes\nabla\varphi])-r(n)=\dot{n}\,,\quad\mbox{ in }\Omega\times[0\,,\tau)\,,\nonumber\\
&D[\nabla n]\mb=0\,,\quad\mbox{ on }\partial\Omega\times[0\,,\tau)\,;\nonumber\\
&-\epsilon\Delta\varphi=\chi|_{{}_{\Omega}}ez\cdot n\,,\quad\mbox{ in }\Real^{3}\times[0\,,\tau)\,,\\
&\ljump\nabla\varphi\rjump\cdot\mb=0\,,\quad\mbox{ on }\partial\Omega\times[0\,,\tau)\,,\nonumber\\
&n_{o}(x)=n(x\,,0)\,,\quad\quad\varphi_{o}(x)=\varphi(x\,,0)\,,\mbox{ in }\Omega\,.\nonumber
\end{eqnarray}

In (\ref{RDD}) the $k\times k$ matrix $N(n)$ is defined as $N(n)=\diag\{n_{1}\,,n_{2}\,,\ldots n_{k}\}$, the $k\times k$ symmetric and semi-definite positive \emph{Diffusion} and \emph{Mobility} matrices $D$ and $M$ are correlated by the Einstein-Smoluchowsky relation:
\begin{equation}\label{einstein}
D=\frac{\theta k_{B}}{e} M\,,
\end{equation}
and $r(n)$ is the $k-$dimensional \emph{recombination} array.

The relation between the matrix $S$ in  (\ref{BVP1}) and the matrix $M$ in (\ref{RDD}) is\footnote{We assume that the mobility is independent on $n$, as pointed out into \cite{KD12}: in such a case $S(n)$ must be restricted to the form
\begin{equation*}
S(n)=S_{o}N(n)\,,
\end{equation*}
where $S_{o}$ is a $k\times k$ matrix with constant components.}
\begin{equation}
M=eS(n)N^{-1}(n)\,,
\end{equation}
where $N^{-1}(n)$ denotes the diagonal $k\times k$ matrix whose entries are $n_{j}^{-1}$ when $n_{j}\neq 0$ and $0$ when $n_{j}=0$. The semi-definite positiveness of $M$ accounts for excitation carriers whose mobilities (and hence by (\ref{einstein}) the associated diffusivities) are either zero or negligible with respect to those of other charge carriers. However in \S.\textbf{3} we shall see that in order to get decay time estimates we shall require the stronger requirement of positive-definiteness for $M$. 

Finally, the $k\times k$ matrix $H(n)$ and the $k-$dimensional array $r(n)$ are related by:
\begin{equation}\label{OnsagerH}
H(n)s=r(n)\,,
\end{equation}
where, by following \cite{MIE11}, we assume that
\begin{equation}\label{Hmatrix}
H(n)=\sum_{h=1}^{s}k_{h}\ell(\frac{ n^{a^{h}} }{ c^{a^{h}}}\,,\frac{n^{b^{h}} }{ c^{b^{h}}})(a^{h}-b^{h})\otimes(a^{h}-b^{h})\,,\quad v^{a}=\prod_{j=1}^{k}v_{j}^{a_{j}}\,.
\end{equation}
In (\ref{Hmatrix}) $h=1\,,2\,,\ldots s$ denote the number of recombination processes with rates $k_{h}$ with the two $k-$dimensional arrays $a^{h}=(a_{1}^{h}\,,a_{2}^{h}\,,\ldots a_{k}^{h})$ and $b^{h}=(b_{1}^{h}\,,b_{2}^{h}\,,\ldots b_{k}^{h})$
describing the $h^{th}$ recombination mechanism with rate $k_{h}$:
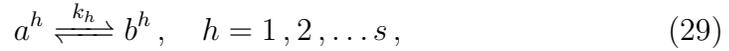
\begin{equation}
a^{h}\overset{k_{h}}{\ce{<=>}}b^{h}\,,\quad h=1\,,2\,,\ldots s\,, 
\end{equation}
and the function $\ell(x\,,y)$ is the logarithmic mean:
\begin{equation}
\ell(x\,,y)=
\begin{cases}
\frac{x-y}{\log x-\log y}\,,\quad x\neq y\,,\\
x\,,\quad x=y\,.
\end{cases}
\end{equation}

Whenever we assume for $F(n)$ the Gibbs entropy (\ref{Gibbs}), then from (\ref{OnsagerH}) and (\ref{Hmatrix}) we arrive at a polynomial expression for the recombination term $r(n)$ (\emph{vid. e.g.} \cite{GMI13} for the details):
\begin{equation}\label{rnbalanced}
r(n)=\sum_{h=1}^{s}k_{h}(\frac{ n^{a^{h}} }{ c^{a^{h}}}-\frac{n^{b^{h}} }{ c^{b^{h}}})(a^{h}-b^{h})\,,
\end{equation}
where we used the identity $\log{v^{a}}=a\cdot\log{v}$. In order to arrive to (\ref{rnbalanced}) we implicitly assumed that all the recombination mechanisms are \emph{detailed balanced} \cite{MIE11}, that is there exists a steady recombination state which in our case coincides with $c$; further, let
\begin{equation}
\caS\equiv\span\{a^{i}-b^{i}\mid i=1,2,\ldots s\}\subset\Real^{k}\,,
\end{equation}
and
\begin{equation}
\caS^{\perp}\equiv\{v\in\Real^{k}\mid v\cdot u\,,\quad \forall u\in\caS\}\,,
\end{equation}
then from (\ref{rnbalanced}) we have
\begin{equation}
r(n)\in\caS\,.
\end{equation}

If we consider now the differential form of the charge conservation (\ref{conservationcharge}), namely:
\begin{equation}\label{conservationcharge1}
\frac{\dif}{\dif t}Q(t)=e\int_{\Omega}z\cdot\dot{n}=0\,,\quad\forall t\in[0\,,\tau)\,,
\end{equation}
then from (\ref{RDD}) we obtain
\begin{equation}
0=\int_{\Omega}r(n)\cdot z=\overline{r(n)}\cdot z\,,
\end{equation}
which in turn, by (\ref{rnbalanced}), implies the electrical neutrality of each recombination mechanism:
\begin{equation}\label{ENRM}
z\cdot(a^{h}-b^{h})=0\,,\quad h=1\,,2\,,\ldots s\,,
\end{equation}
and hence $z\in\caS^{\perp}$ in such a way that the local orthogonality condition holds:
\begin{equation}\label{firstinte}
r(n)\cdot z=0\,.
\end{equation}

In \cite{MIE11} it is remarked that the term $r(n)$ can be represented as a polynomial relation in $n$ only if we identify $F(n)$ with the Gibbs entropy (\ref{Gibbs}), whereas for a different entropic term like \emph{e.g.\/} the Fermi-Dirac potentials this is not possible.

Equation (\ref{RDD})$_{1}$, which represents the conservation of electric current normalized with respect to $e$,  was first proposed for scintillators in \cite{VA08} by following \cite{FO64} and \cite{AR66} and with $r(n)$ a polynomial, at most cubic, function of $n$; it was used for scintillators into \cite{LGW11},  \cite{WGLU11}, \cite{KD12}-\cite{SW15a}. Moreover, special cases of this equations were widely used to model scintillation at a phenomenological level in many theoretical and experimental paper, as we shall describe in details in \S.\textbf{4}. We also remark that the boundary value problem (\ref{RDD}) is the same obtained, by starting from a different approach and with a different reaction term $r(n)$, in \cite{AGH02} for semiconductors (\emph{vid.\/ also} \cite{MIE11}, \cite{MIE15}).

In the available phenomenological models for scintillation the recombination term $r(n)$ is generally assumed as a cubic expression in $n$:
\begin{equation}\label{recombination}
r_{i}(n)=r_{i}^{0}+\sum_{j=1}^{k}A_{ij}n_{j}+\sum_{h,j=1}^{k}B_{ijh}n_{h}n_{j}+\sum_{h,m,j=1}^{k}C_{ijhm}n_{h}n_{m}n_{j}\,,\quad i=1\,,2\,,\ldots\,,k\,.
\end{equation}
The terms $r_{i}^{o}$ describes the excitation carriers creation rate in the crystal under irradiation; the other terms in (\ref{recombination}) are further splitted in order to represent different recombination mechanism. The terms $A_{ij}$, which accounts for linear recombination, can be further decomposed into three terms
\begin{equation}\label{linearcomb}
A_{ij}=A^{r}_{ij}+A^{nr}_{ij}+A^{e}_{ij}\,,
\end{equation}
which represent respectively radiative recombination with photon emission, non-radiative recombination without photon emission and exchange between excitation carriers. The quadratic recombination is represented by the terms $B_{ijh}$ which can also be splitted into $B_{ijh}=B^{r}_{ijh}+B^{nr}_{ijh}$ with the same meaning of superscript as in (\ref{linearcomb}), whereas the third-order or Auger totally non-radiative recombination is described by the terms $C_{ijhm}=C^{nr}_{ijhm}$.

In order to reconcile (\ref{recombination}) with (\ref{rnbalanced}) we consider here the following example for $k=3$ where $n_{1}=n_{e}$ and $n_{2}=n_{h}$ represent respectively the electrons and holes densities and $n_{3}=n_{ex}$ represents the excitons density. By following, for instance, the description provided in  \S.3.3 of \cite{KTV20} then we may have these recombination mechanisms:
\begin{enumerate}
\item $n_{1}+n_{2}=0\,,\quad a^{1}=(1\,,1\,,0)\,,\quad b^{1}=(0\,,0\,,0)$\,,
\item $n_{1}+n_{2}=n_{3}\,,\quad a^{2}=(1\,,1\,,0)\,,\quad b^{2}=(0\,,0\,,1)$\,,
\item $n_{1}+n_{2}=n_{1}+n_{2}+n_{3}\,,\quad a^{3}=(1\,,1\,,0)\,,\quad b^{3}=(1\,,1\,,1)$\,,
\item $n_{1}=2n_{1}+n_{2}\,,\quad a^{4}=(1\,,0\,,0)\,,\quad b^{4}=(2\,,1\,,0)$\,,
\item $n_{1}=n_{1}+n_{3}\,,\quad a^{5}=(1\,,0\,,0)\,,\quad b^{5}=(1\,,0\,,1)$\,,
\item $n_{2}=n_{1}+2n_{2}\,,\quad a^{6}=(0\,,1\,,0)\,,\quad b^{6}=(1\,,2\,,0)$\,,
\item $n_{2}=n_{2}+n_{3}\,,\quad a^{7}=(0\,,1\,,0)\,,\quad b^{7}=(0\,,1\,,1)$\,,
\end{enumerate}
and then from (\ref{rnbalanced}) we get:
\begin{eqnarray}
r_{1}(n)=r_{2}(n)&=&-k_{1}-\frac{k_{4}}{c_{1}}n_{1}-\frac{k_{6}}{c_{2}}n_{2}-\frac{k_{2}}{c_{3}}n_{3}\nonumber\\
&+&\frac{k_{1}+k_{2}}{c_{1}c_{2}}n_{1}n_{2}+\frac{k_{4}}{c^{2}_{1}c_{2}}n^{2}_{1}n_{2}+\frac{k_{6}}{c_{1}c^{2}_{2}}n_{1}n^{2}_{2}\,,\\
r_{3}(n)&=&-\frac{k_{3}}{c_{1}}n_{1}-\frac{k_{7}}{c_{2}}n_{2}+\frac{k_{2}}{c_{3}}n_{3}\nonumber\\
&+&\frac{k_{2}-k_{3}}{c_{1}c_{2}}n_{1}n_{2}+\frac{k_{5}}{c_{1}c_{3}}n_{1}n_{3}+\frac{k_{7}}{c_{2}c_{3}}n_{2}n_{3}\,,\nonumber
\end{eqnarray}
which can be trivially put into the form (\ref{recombination}), with the appropriate identification of radiative, non-radiative or exchange terms: we notice that the mechanisms 3, 4 and 6 are Auger recombination, whereas 2, 5 and 7 represents scattering and the mechanism 1 is the simple electron-hole recombination  (\emph{cf. e.g.\/} the models proposed either into \cite{BM09} or \cite{BMSVW09} where $n_{e}=n_{h}=n_{eh}$).

As we already remarked more complex expressions for (\ref{recombination}) can be proposed: for instance into \cite{GBD15}, for $k=7$, we had:
\begin{eqnarray}
r_{1}(n)&=&r_{1}^{0}+A_{14}n_{4}+B_{113}n_{1}n_{3}+B_{115}n_{1}n_{5}+B_{117}n_{1}n_{7}\,,\nonumber\\
r_{2}(n)&=&r_{2}^{0}+A_{23}n_{3}+A_{25}n_{5}+B_{223}n_{2}n_{3}+B_{224}n_{2}n_{4}+B_{227}n_{2}n_{7}\,,\nonumber\\
r_{3}(n)&=&A_{33}n_{3}+B_{323}n_{2}n_{3}+B_{336}n_{3}n_{6}\,,\nonumber\\
r_{4}(n)&=&A_{44}n_{4}+B_{417}n_{1}n_{7}+B_{424}n_{2}n_{4}\,,\\
r_{5}(n)&=&A_{55}n_{5}+B_{515}n_{1}n_{5}+B_{527}n_{2}n_{7}\,,\nonumber\\
r_{6}(n)&=&r_{6}^{0}+A_{66}n_{6}+B_{613}n_{1}n_{3}\,,\nonumber\\
r_{7}(n)&=&A_{77}n_{7}+B_{715}n_{1}n_{5}+B_{724}n_{2}n_{4}\,,\nonumber
\end{eqnarray}
with $n_{1}$ and $n_{2}$ the  electron and holes densities, $n_{3}$ and $n_{6}$ the self-trapped electrons and excitons respectively, $n_{4}$ and $n_{5}$ the electron and holes captured by the activation centers and finally $n_{7}$ denotes the excited activator centers: in \cite{GBD15} a detailed description of the physical motivation of these relations is provided. We notice that in this model the Auger mechanism is missing and the recombination terms are at most quadratic.

\subsection{Stationary solutions}

It is trivial to show that the stationary solutions, that is equilibrium solutions $n^{\infty}$ for the boundary value problem (\ref{RDD}) with $\dot{n}=0$, can be obtained by setting $s=0$, which by (\ref{OnsagerH}) leads to the equilibrium condition:
\begin{equation}\label{stationaryreaction}
r(n^{\infty})=0\,.
\end{equation}
When $F(n)$ is identified with the Gibbs entropy (\ref{Gibbs}), then from $s=0$ we have
\begin{equation}\label{stationary1}
s_{j}=ez_{j}\varphi^{\infty}+\theta k_{B}\log\frac{n_{j}^{\infty}}{c_{j}}=0\,,\quad j=1\,,2\,,\ldots\,,k\,,
\end{equation}
where the stationary electric field $\varphi^{\infty}$ is the solution of:
\begin{equation}\label{laplaceinfty}
-\epsilon\Delta\varphi^{\infty}=e\sum_{j=1}^{k}z_{j}n_{j}^{\infty}\,,\mbox{ in }\Omega\,,
\end{equation}
with Neumann boundary conditions on $\partial\Omega$. From (\ref{stationary1}) then we have
\begin{equation}\label{stationary2}
n_{j}^{\infty}(x)=c_{j}\exp(-\frac{ez_{j}\varphi^{\infty}(x)}{\theta k_{B}})\,,\quad j=1\,,2\,,\ldots\,,k\,,
\end{equation}
and (\ref{laplaceinfty}), (\ref{stationary2}) together leads to the semilinear Poisson-Boltzmann equation for the stationary electric field:
\begin{equation}\label{LBG}
-\epsilon\Delta\varphi^{\infty}=e\sum_{j=1}^{k}z_{j}c_{j}\exp(-\frac{ez_{j}\varphi^{\infty}}{\theta k_{B}})\,,\mbox{ in }\Omega\,;
\end{equation}
into \cite{LI09} an uniqueness result in $H^{1}(\Omega)$ for (\ref{LBG}) with Dirichlet boundary conditions was obtained.

It is trivial to show that the equilibrium solution $n^{\infty}$ is also a steady state $c$ for the detailed balance condition and indeed in (\ref{OnsagerH}) we have:
\begin{equation}
\frac{{(n^{\infty})}^{a^{h}}}{c^{a^{h}}}-\frac{{(n^{\infty})}^{b^{h}}}{c^{b^{h}}}=\exp(-\frac{e\varphi}{k_{B}\theta}(z\cdot a^{h}-z\cdot b^{h}))=0\,,
\end{equation}
by (\ref{stationary2}) and (\ref{ENRM}), and hence the equilibrium condition (\ref{stationaryreaction}) is satisfied.

We finally deal with the problem of the determination of $c=(c_{1}\,,c_{2}\,,\ldots\,,c_{k})$ in (\ref{stationary2}): if we define the $k\times k$ diagonal matrix $L(x)$
\begin{equation}
L(x)=\diag\{\exp(-\frac{ez_{1}\varphi^{\infty}(x)}{\theta k_{B}})\,,\exp(-\frac{ez_{2}\varphi^{\infty}(x)}{\theta k_{B}})\,,\ldots,\exp(-\frac{ez_{k}\varphi^{\infty}(x)}{\theta k_{B}})\}\,,
\end{equation}
then (\ref{stationary2}) can be written as:
\begin{equation}\label{LC}
n^{\infty}(x)=L(x)c\,.
\end{equation}
and then, by (\ref{conservationcharge}) $\overline{L(x)}c\in\caS$. 

To obtain the explicit value of $c$ we may use (\ref{stationaryreaction}), whereas the uniqueness of $c$ follows instead, as pointed out in \cite{FK18}, by monotonicity, (\ref{conservationcharge}) and the uniqueness  of $\varphi^{\infty}\in H^{1}(\Omega)$ with $\bar{\varphi}^{\infty}=0$.

\section{Existence and asymptotic decay}

In this section we shall give an account of the existing results concerning the existence of solutions and the asymptotic estimates for the boundary value problem (\ref{RDD}): we remark that the latter results are important in order to get a scintillation decay time estimate. We shall not enter into the mathematical details which can be found in the references we quote, rather we shall adapt if necessary these results to the specific cases of our boundary value problems. 

The problem of finding existence, asymptotic estimates and qualitative bounds for the solutions for the coupled boundary value problem (\ref{RDD}) has received a strong attention in the recent years, \emph{vid. e.g.\/} \cite{CHJU17}, \cite{CHJU18},\cite{HHMM18},  \cite{GA94}-\cite{FK20} and the many references quoted therein: to this regard it is important to remark that most of these results deal with semiconductors or chemical reactions which differ from scintillators by the structure of the reaction term $r(n)$. Most of these results are based on the so-called \emph{Entropy method},  whose importance is explained in full in these words taken from \cite{DFM08}:
\begin{quote}
\emph{The entropy method refers to the general idea of a functional inequality relationship between an entropy functional of a system and its monotone change in time, usually called the entropy dissipation. Such an entropy-entropy dissipation inequality entails convergence to an entropy minimizing equilibrium state, at first in entropy and further in $L^{1}$ using Czisz\'ar-Kullback-Pinsker-type inequalities. The entropy approach is per se a nonlinear method avoiding any kind of linearization and capable of providing explicitly computable convergence rates. Moreover, being based on functional inequalities rather than particular differential equations, it has the advantage of being quite robust with respect to model variations.}
\end{quote}

In the next subsection we shall show how some of these results can be extended to the RDD equations for scintillators.

As far as the decay time is concerned, the available experimental data (\emph{vid. e.g.} the recent analysis in \cite{SWI14}) and the numerical solution of phenomenological models as in \cite{LGW17}, show that the excitation carriers decay exponentially in time to an asymptotic value $n_{\infty}$, namely:
\begin{equation}
\|n(\cdot\,,t)-n_{\infty}(\cdot)\|=A_{f}\exp(-t/\tau_{f})+A_{s}\exp(-t/\tau_{s})\,,
\end{equation}
where the indeces $f$ and $s$ denotes the so-called \emph{fast} and \emph{slow} components of the excitation, respectively. Accordingly, since by definition the \emph{Decay time} is the time required for scintillation emission to decrease to $e^{-1}$ of its maximum, then we get a Fast Decay Time $\tau_{f}$ and a Slow Decay Time $\tau_{s}$. In many cases one of the components is negligible and the decay obeys a simple exponential law, which can be also used to describe an average decay time. 

\subsection{Global existence}

The first results concerning existence theory for the boundary value problems like (\ref{RDD}) was obtained in \cite{FI15}, \cite{FI17}: such a result, which relies on the notion of global \emph{renormalised solutions} leaves still open, as pointed out in \cite{FI17}, the problem the existence  of weak or even smooth global solutions in time: the reason, as pointed out in detail into \cite{FI17} are the growth condition on the reaction/recombination terms. 

These results were the extended into \cite{CHJU17} and \cite{CHJU18} and here we shall follow the latter. The main strongpoint of these results is that no growth condition are imposed on the reaction/recombination term. In our case the diffusion is linear: unfortunately we cannot extend their results to the non-linear case they study, because for (\ref{einstein}) this would imply also a non linear mobility, whereas in all these paper the mobility is implicitly assumed $M=I$ with $I$ the $k\times k$ identity.

First of all we recall the notion of renormalised solution, first introduced into \cite{DL88}-\cite{DL89b} for the Boltzmann and transport equations, as it was given into \cite{FI15}. 

An excitation density vector $n$ is a renormalised solutions for (\ref{RDD}) if for all functions $\xi:\caM\rightarrow\Real$  with compactly supported derivative $\nabla_{n}\xi$ , the function $\xi(n)$ must satisfy the equation derived from (\ref{RDD}) by a formal application of the chain rule in a weak sense. As it is pointed out into \cite{VIL02}, the function $\xi$ must belongs to a well-choosen class of admissible solutions. The  physical interpretation of these renormalised solution is that they gives a distributional sense to the boundary value problem (\ref{BVPadim}).

More precisely we say that:
\begin{itemize}
\item $n=(n_{1}\,,n_{2}\,,\ldots\,,n_{k})$  is a \emph{renormalized solutions} for (\ref{RDD}) if $\forall \tau>0$, $n_{i}\in L^{2}(H^{1}(\Omega);[0\,,\tau))$ and for any $\xi\in C^{\infty}(\caM)$, such that $\nabla_{n}\xi\in C_{0}^{\infty}(\caM; \Real^{k})$ and $\psi\in C_{0}^{\infty}(\bar{\Omega}\times [0\,,\tau))$, it holds:
\begin{eqnarray}\label{reno1}
&&\int_{0}^{\tau}\int_{\Omega}\xi(n)\dot{\psi}=\nonumber\\
&=&\int_{0}^{\tau}\int_{\Omega}([D\nabla n+MN(n)z\otimes\nabla\varphi]\cdot\nabla_{n}\nabla_{n}\xi[\nabla n]+r(n)\cdot\nabla_{n}\xi)\psi\nonumber\\
&+&\int_{0}^{\tau}\int_{\Omega}(D\nabla n+MN(n)z\otimes\nabla\varphi)\cdot\nabla_{n}\xi\otimes\nabla\psi\,.
\end{eqnarray}
\end{itemize}

We leave out all the details and recall only the main hypotheses and results given in \cite{CHJU18} and first of all we assume that there exist numbers $\pi_{i}>0$ and $\lambda_{i}\in\Real$, $i=1,2,\ldots,k$ such that $\forall n\equiv(n_{1}\,,n_{2}\,,\ldots\,,n_{k})\in (0\,,\infty)^{k}$, the following inequality holds:\footnote{A certain care is requested when we look at (\ref{entropyineq}) and (\ref{quasipos}), since the reaction term $f(u)$ in \cite{CHJU18} is the opposite of our recombination term $r(n)$, say $f(u)=-r(n)$. Hence the reversed inequalities and different definition than those given into \cite{CHJU18}. }
\begin{equation}\label{entropyineq}
\sum_{i=1}^{k}\pi_{i}r_{i}(n)(\log\frac{n_{i}}{c_{i}}+\lambda_{i})\geq 0\,,
\end{equation}
a condition which implies the \emph{quasi-negativity} of $r(n)$, that is (\emph{cf.} the models in \cite{GBD15} and \cite{BM09}):
\begin{equation}\label{quasipos}
r_{i}(n_{1}\,,\ldots, n_{i-1}\,,0\,,n_{i+1}\,,\ldots,n_{k})\leq 0\,,\quad\forall n\in\caM\,,\forall i=1\,,2\,,\ldots, k\,,
\end{equation}
which grants the non-negativity of solutions.

Condition (\ref{entropyineq}) further ensure the existence of a so-called \emph{total entropy}:
\begin{equation}\label{entropytotal}
\caH(n)=\int_{\Omega}\sum_{i=1}^{k}\pi_{i}n_{i}(\log\frac{n_{i}}{c_{i}}-1+\lambda_{i})+\exp(-\lambda_{i})\,,
\end{equation}
which is a Lyapunov functional for the reaction system (\ref{kinetic1})  if $\pi_{i}=1\,,\forall i$.

Provided these preliminary conditions are satisfied, then the main hypotheses from \cite{CHJU18}  can be rephrased, within the context of our treatment, as:
\begin{itemize}
\item [(H1)] Drift term: $\nabla\varphi\in L^{\infty}([0\,,\tau)\,; L^{\infty}(\Omega\,,\Real^{k+3}))$;
\item [(H2)]
\begin{itemize}
\item[i)] Recombination term: $r(n):\caM\rightarrow\Real^{k}$ is locally Lipschitz continuos, that is there exists a  function $K(\cdot):[0\,,+\infty)\rightarrow [0\,,+\infty)$ non-decreasing and such that a.e. $(x\,,t)\in\Omega\times[0\,,\tau)$ and $\forall n\,,\hat{n}\in\caM$:
\begin{equation}\label{lips}
\|r(n)-r(\hat{n})\|\leq K(\max\{\|n\|\,,\|\hat{n}\|\})\|n-\hat{n}\|\,,
\end{equation}
and $\forall\hat{\tau}>0$, then $r(0)\in L^{2}(\Omega\times[0\,,\hat{\tau}))$.
\item[ii)] the inequality (\ref{entropyineq}) holds\,,
\item[iii)] there exist $m\in\Natural$ such that $\forall n\in\caM$ with $\sum_{j=1}^{k}n_{j}\geq m$, then $\sum_{j=1}^{k}r_{j}(n)\geq 0$\,;
\end{itemize}
\item [(H3)] Initial data: $n_{o}=(n^{0}_{1}\,,n^{0}_{2}\,,\ldots\,,n^{0}_{k})\in L^{\infty}(\Omega\,,\Real^{k})$, such that $\inf n^{o}_{j}>0\,,j=1\,,2\,,\ldots, k$;
\item [(H4)] The mobility matrix $M$ (and hence for (\ref{einstein}) the diffusion matrix $D$) is diagonal and positive-definite.
\end{itemize}

From a physical point of view, the first three hypotheses requires simply a certain degree of regularity on the initial data, the electric field and the recombination term, whereas the last one rules out bot cross-mobility (and diffusion), as well as the possibility to deal with charge carriers with no mobility.

A first consequence of the hypotheses (H1)-(H4) is that there exists
\begin{itemize}
\item a global renormalised solution $n=(n_{1}\,,n_{2}\,,\ldots\,,n_{k})$ satisfying $n_{i}\geq 0$ in $\Omega$, $i=1,2,\ldots k$\,, $\forall t$\,,
\item $\caH(n)=<+\infty$\,, $\forall t$\,;
\end{itemize}
moreover they imply the main result of \cite{CHJU18} (\emph{Weak-strong uniqueness of the solutions}) which states that the boundary value problem (\ref{RDD}):
\begin{itemize}
\item admits a weak-strong renormalised solution $n=(n_{1}\,,n_{2}\,,\ldots\,,n_{k})$ satisfying $n_{i}\geq 0$ in $\Omega\times[0\,,\tau)$, $i=1,2,\ldots k$\,, $\forall t$\,,
\item $\caH(n)=<+\infty$\,, in $\Omega\times[0\,,\tau)$\,,
\item $n\in (0\,,\tau\,;L^{1}(\Omega))$, \emph{cf.\/} (\ref{L1norm})\,;
\end{itemize}
Finally, for $n$ a renormalized solution to (\ref{BVPadim}) and $v=(v_{1}\,,\ldots v_{k})$ is a \emph{strong} solution to (\ref{BVPadim}) on some time interval $[0\,,\tau^{\star})$ with $\tau^{\star}\leq\tau$ in the following sense: there exist $C>c>0$ such that:
\begin{itemize}
\item $c\leq v_{i}(x\,,t)\leq C\,,\quad(x\,,t)\in\Omega\times[0\,,\tau^{\star})\,, i=1\,,2\,,\ldots,k$;
\item $\|\dot{v}\|_{L^{\infty}(\Omega\times[0\,,\tau^{\star}))}+\|\nabla v\|_{L^{\infty}(\Omega\times[0\,,\tau^{\star}))}\leq C$\,,
\end{itemize}
then if, for any $s\in(0\,,\tau^{\star})$,  $\phi\in C^{\infty}(\overline{\Omega}\times[0\,,\tau^{\star}))$:
\begin{equation}\label{weakRDD}
\int_{0}^{s}\int_{\Omega}\phi\dot{v}=-\int_{0}^{s}\int_{\Omega}(D[\nabla v]+MN(v)z\otimes\nabla\varphi)\cdot\nabla\phi-\int_{0}^{s}\int_{\Omega}\phi r(v)\,;
\end{equation}
then $n(s\,,t)=v(s\,,t)$, for $x\in\Omega$, $s\in (0\,,\tau^{\star})$.

We remember that, as pointed out into \cite{FI17}, we can say nothing about the global existence in time of smooth solutions: however this result is important in order to get a general framework for numerical solutions of (\ref{BVPadim}).

\subsection{Asymptotic estimate for the decay time}

In \cite{FK18} an explicit estimate of the asymptotic convergence was obtained for the Rosbroeck model for semiconductors with Shockley-Read-Hall potential and $k=2$: here we shall show how the results obtained there can be adapted to the case of scintillators in order to obtain an explicit estimate for the decay time. Once again we shall leave out the technicalities and details and we shall give only the main results using our language and notations.

We have already seen that we can safely assume uniqueness for $c$: given this in \cite{FK18} some preliminary bounds are necessary, namely the following for $c=(c_{1}\,,c_{2})$ and $n^{\infty}=(n_{1}^{\infty}\,,n_{2}^{\infty})$:
\begin{equation}\label{bounds}	
c_{j}\leq e^{\Phi^{\infty}}\,,\quad n^{\infty}_{j}\leq e^{2\Phi^{\infty}}\,,\quad (n^{\infty}_{j})^{-1}\leq e^{2\Phi^{\infty}}\,,\quad j=1\,,2\,,
\end{equation}
with
\begin{equation}
\Phi^{\infty}=\|\frac{ez\varphi^{\infty}}{\theta k_{B}}\|_{L^{\infty}(\Omega)}\,.
\end{equation}

Moreover, for the recombination $r(n)$ given by (\ref{recombination}) there exists a constant $k_{o}$ such that:
\begin{equation}\label{bound0}
0<k_{o}=\|r(0)\|_{L^{\infty}(\Omega)}\leq\|r(n)\|_{L^{\infty}(\Omega)}\,.
\end{equation}
If we consider the free-energy $U(n)$ defined by (\ref{electroenergy}) with the choice of the Gibbs entropy (\ref{Gibbs}) for the term $F(n)$, then the main results of \cite{FK18} are based on the derived notion of \emph{Relative Gibbs free-energy} (or \emph{relative entropy} according to \cite{FK18}): 
\begin{equation}
\caG(u\mid v)=U(u)-U(v)-\caD U(v)(u-v)\,;
\end{equation}
by an explicit calculation it can be show that
\begin{equation}\label{relativeentropy1}
\caG(n\mid n^{\infty})=\int_{\Omega}\sum_{i=1}^{k}n_{i}\log(\frac{n_{i}}{n^{\infty}_{i}})-(n_{i}-n^{\infty}_{i})+\frac{1}{2}\epsilon\|\nabla\varphi-\nabla\varphi^{\infty}\|^{2}\,,
\end{equation}
and then, by an easy calculation, it can be shown that the Dissipation $\caD$ defined by (\ref{dissipation1}) is given by 
\begin{equation}\label{dissipation}
\Diss(n\,,\varphi)=-\frac{\dif}{\dif t}\caG(n\mid n^{\infty})\,.
\end{equation}

In \cite{FK18} by starting from (\ref{dissipation}) and by the means of a repeated use of Csisz\'ar-Kullback-Pinsker type inequalities, provided (\ref{bounds}) hold, the two following estimates were obtained:
\begin{eqnarray}\label{estimate0}
&\Diss(n\,,\varphi)\geq C_{1}\,\caG(n\,,\varphi)\,,\nonumber\\
\\
&\|n-n_{\infty}\|^{2}_{L^{1}(\Omega)}+\|\varphi-\varphi_{\infty}\|^{2}_{H^{1}(\Omega)}\leq C_{2}\,\caG(n_{o}\,,\varphi_{o})\,\exp(-C_{1}t)\,,\nonumber
\end{eqnarray}
where $\varphi_{o}$ is the unique solution of (\ref{laplace}) for the initial data $n_{o}$.

The most remarkable feature of this results is that both $C_{1}$ and $C_{2}$ have an explicit dependence on the parameters of (\ref{BVPadim}):
\begin{eqnarray}\label{parameter}
C_{1}^{-1}&=&\frac{1}{2}\exp({2\Phi^{\infty}})\max\{\frac{T}{m}\exp({2\Phi^{\infty}})\,,\frac{1}{k_{o}}\}\cdot(1+\caL(\Omega)\exp({2\Phi^{\infty}}))\,,\nonumber\\
C_{2}&=&3\exp({2\Phi^{\infty}})+\frac{1}{2}\caG(n_{o}\,,\varphi_{o})+2(1+\caL(\Omega))\,;
\end{eqnarray}
in (\ref{parameter}) $\caL(\Omega)>0$ is the constant in the Poincar\'e inequality:
\begin{equation}
\|\Psi\|^{2}_{L^{p}(\Omega)}\leq\caL(\Omega)\|\nabla\Psi\|^{2}_{L^{p}(\Omega)}\,,\quad\forall\Psi\in H^{1}(\Omega)\,,\bar{\Psi}=0\,,
\end{equation}
which in $L^{2}(\Omega)$ reduces to $\caL(\Omega)=\lambda_{1}^{-1}$ with $\lambda_{1}$ the first eigenvalue of:
\begin{eqnarray}
\Delta\Psi=\lambda\Psi\,,\quad\mbox{ in }\Omega\,,\\
\Psi=0\,,\quad\mbox{ on }\partial\Omega\,.\nonumber
\end{eqnarray}
We notice that, for $\Omega$ a sphere of radius $R$ it is $\caL(\Omega)\leq(2R/\pi)^{2}$ \cite{PW60}: for $L=2R$ then we have $\caL(\Omega)\leq 0.1$. Such a value is consistent with the value $\caL(\Omega)\approx 0.07$ which was obtained into \cite{D3} where the $\lambda_{1}$ was calculated as the the square root of the reciprocal of the first zero of the Bessel function $J'_{o}$.

The expression for the decay time $\tau=C_{1}^{-1}$ depends, by (\ref{parameter})$_{1}$, in an explicit manner on the mobility parameter $m$, the reaction term $r(0)$, the initial data $n_{o}$,  and the scintillation volume $\Omega$: as far as we know it is the first explicit estimate of scintillator decay time which depend on the (measurable) constitutive parameters of the model, albeit limited to the case $k=2$. In \cite{D3} we showed how, for four different scintillators these results gave a very good estimate for the experimentally measured fast decay time.

\section{The approximated models}

\subsection{The adimensionalised boundary value problem}

In order to put the boundary value problem (\ref{RDD}) in an adimensionalised form, we begin with the choice of a characteristic length and time pair $(L\,,T)$ related to the space and time scintillation scales, in order to define the dimensionless coordinates $(\zeta\,,\tau)$:
\begin{equation}
\zeta=\frac{x}{L}\,,\quad \tau=\frac{t}{T}\,,
\end{equation}
and the dimensionless excitation carriers density $u=u(\zeta\,,\tau)$ and electric potential $\psi=\psi(\zeta\,,\tau)$:
\begin{equation}
u=nL^{3}\,,\quad\psi=\varphi\frac{\epsilon}{eL}\,.
\end{equation}
If we further set:
\begin{equation}
M^{*}=\mu^{-1}M\,,\quad r^{*}=\kappa^{-1}r\,,
\end{equation}
where $\mu$ is the greatest eigenvalues of $M$ and and $\kappa$ the greater component of $r$ or, as in \cite{D3}, the value of their norm, then equations (\ref{RDD})$_{1,3}$ can be rendered dimensionless
\begin{eqnarray}\label{BVPadim}
&\dive{}_{\zeta}(dD^{*}[\nabla_{\zeta}u]+mM^{*}[ez\otimes\nabla_{\zeta}\psi])-kr^{*}(u)=u_{\tau}\,,\\
&-\Delta_{\zeta}\psi=z\cdot u\,,\nonumber
\end{eqnarray}
and dependent on the three adimensional parameters:
\begin{equation}
d=\frac{Tk_{B}\theta}{eL^{2}}\mu\,,\quad m=\frac{eT}{\epsilon L^{3}}\mu\,,\quad k=\kappa T\,.
\end{equation}

We remark that the parameter $k$ depends on the incoming energy by the means of $r(n)$, whereas $m$ depends only on constituive or scaling parameters: accordingly, depending on the energy of ionizing radiation we may have different physically meaningful regimes which are described in the following subsections. Further, for 
\begin{equation}
\delta=\frac{k_{B}\theta}{e}\mu\,,
\end{equation}
the greatest eigenvalue of $D$ and provided we define the \emph{diffusion length} $L_{D}=\sqrt{\delta\,T}$, then we have
\begin{equation}\label{diffL}
d=(\frac{L_{D}}{L})^{2}\,;
\end{equation}
we notice that an experimental estimate of $L_{D}$ is given \emph{e.g.\/} into \cite{BM09}.

The set of dimensionless parameters $\{d\,,m\,,k\}$ which appears into (\ref{BVPadim})$_{1}$ describes three main regimes for the boundary value problem, depending on the crystal constitutive properties and on the initial ionizing energy. In the ext subsections we shall show how many of the most used phenomenological models for scintillators, can be encompassed within (\ref{BVPadim})$_{1}$ by an appropriate choice of these parameters. For these models, which we define approximated because some terms of (\ref{BVPadim})$_{1}$ can be neglected, we shall also briefly describe the results concerning existence and decay time which are available in the literature.

\subsection{The Reaction-Diffusion approximation}

Let $d\approx k=O(1)$, that is:
\begin{equation}
\kappa\approx\frac{\theta k_{B}}{eL^{2}}\mu=O(1)\,,
\end{equation}
and let
\begin{equation}
\frac{m}{k}=\frac{e}{\epsilon L^{3}}\frac{\mu}{\kappa}=o(\varepsilon)\,,\quad\frac{m}{d}=\frac{e^{2}}{\epsilon L\theta k_{B}}=o(\varepsilon)\,,
\end{equation}
then from (\ref{BVPadim})$_{1}$ we recover, to within higher order terms, the reaction-diffusion equation
\begin{equation}\label{ReaDiff}
\dive D[\nabla n]-r(n)=\dot{n}\,,
\end{equation}
and (\ref{BVPadim})$_{2}$ is uncoupled, the excitation density vector being a data; this phenomenological model, which traces back its ancestry to the analogous models for chemical reactions \cite{KOKU88}, \cite{KOKU96}, is used in many papers dealing with experimental identification of scintillator properties, as in \cite{WLGUF11}-\cite{CZZ20} and did not take into account the drift contribution induced by the electric field $\varphi$.

Reaction-Diffusion equations were studied with a great interests starting from the kinetics of chemical reactions \cite{KOKU96} and there are treatises and textbook devoted to various aspects of them like, \emph{e.g.\/}, \cite{NI11}-\cite{RO84} and many others. However, problems like the existence of solutions and their asymptotic decay have attracted a growing number of studies in recent years, mainly because the difficulties posed by the the lack of control of the recombination terms as pointed out in details into \cite{PI10}-\cite{FI17}.  Most of  the recent results dealing with asymptotic decay deal indeed with reaction terms with quadratic growth for chemical reactions \emph{vid.\/} \cite{DEFE06}-\cite{FETA17} and also, for a different point of view, \cite{MIE11}, \cite{DT18}-\cite{HHMM18}. A first result concerning the existence of renormalised solutions was presented into \cite{DFPV07}, whereas the most general result at the present available is \cite{FLQ19}, where however the nonlinear case of diffusion in porous media is treated. See also the references into \cite{FHKM20}, where these results were extended to the non-isothermal case.

Here we show how these results can be used in the context of our problem, described by the equation (\ref{ReaDiff}). There are, at the best of our knowledge, no results for a general recombination term: however, if we neglect the cubic Auger effect, we can use the existence results of \cite{DFPV07} and decay estimates obtained into \cite{FETA17}.

We shall give only a brief survey of both these results and how their hypotheses fit within the physics underlying scintillation: we leave out all the technical details contained in the cited works. We remark that a complete existence and decay estimate results for the scintillation models has yet to be done and the task, given the nature of the reactive term, is a far from an easy one. A further result presented into \cite{FETA17} is concerned with the existence of renormalized solutions for (\ref{ReaDiff}), whose structure is the same as in (\ref{reno1}) when we neglect the drift terms.

\subsubsection{Existence of global weak solutions}

The main result of \cite{DFPV07} concerns the existence of \emph{global weak solutions} for a reaction/recombination term which is at most quadratic in $n$. The relevant hyoptheses are that the recombination term is locally Lipschitz continuous and quasi-negative (the hypotheses (H2)-i) and (\ref{quasipos}) of \S.\textbf{3}): then, provided $D=\diag\{D_{1}\,,D_{2}\,,\ldots, D_{k}\}$ with $D_{j}>D_{0}>0$, $j=1\,,2\,,\ldots k$ and provided there exists a function $\Theta(x\,,t)\in L^{2}(\Omega\times[0\,,\tau))$, a scalar $\mu\in(0\,,+\infty)$ and a function $K(\cdot)$ defined as in (\ref{lips}) such that:
\begin{itemize}
\item $K(\|r\|)\leq C(1+\|r\|)$\,,
\item $\forall n\in(1\,,+\infty)^{k}$ and a.e. $(x\,,t)$:
\begin{equation}
-\sum_{j=1}^{k}\log(n_{j})r_{j}(n)\leq\Theta+\mu\sum_{j=1}^{k}n_{j}\,,
\end{equation}
\end{itemize}
then (\ref{ReaDiff}) admits a global weak solution for any non-negative initial data such that:
\begin{equation}
\|n_{o}\|\log(\|n_{o}\|)\in L^{2}(\Omega)\,.
\end{equation}

\subsubsection{Asymptotic decay estimates}

The result obtained into \cite{FETA17} is tailored on chemical reactions, which implies mass conservation and detailed or complex balance of reactive terms, conditions which have no correspondence in the physics of scintillation. In this paper, by using arguments and tools which are unsurprisingly related to those used into \cite{FK18}, it is shown that, for $n^{\infty}$ a detailed balanced equilibrium solution corresponding to $r(n^{\infty})=0$, then:
\begin{equation}
\|n-n^{\infty}\|^{2}_{L^{1}(\Omega)}\leq C_{1}\caG(n_{o}\mid n^{\infty})\exp(-C_{2} t)\,,
\end{equation}
where the relative entropy $\caG$ is defined as (\emph{cf.\/} (\ref{relativeentropy1})):
\begin{equation}\label{relativeentropy2}
\caG(n_{o}\mid n^{\infty})=\int_{\Omega}\sum_{i=1}^{k}n_{oi}\log(\frac{n_{oi}}{n^{\infty}_{i}})-(n_{oi}-n^{\infty}_{i})\,,
\end{equation}
$C_{1}=C_{CKP}^{-1}$ is the constant in a Csisz\'ar-Kullback-Pinsker type inequality, whereas the constant $C_{2}$, whose inverse is the estimate for the decay time, is given by:
\begin{equation}\label{C2RD}
C_{2}=\frac{1}{2}\min\{\lambda_{1}\,,\frac{K_{2}H_{1}(t)}{K_{1}}\}\,;
\end{equation}
in this relation $\lambda_{1}=C_{LSI}\min\{D_{i}\}$, with $D=\diag\{D_{1}\,,D_{2}\,,\ldots D_{k}\}$ and $C_{LSI}$ is the constant in the logarithmic Sobolev inequality, whereas the two constants $K_{1,2}$ depends explicitly on $\Omega$, $D$, the set $k_{h}\,,a^{h}\,,b^{h}$, $n^{\infty}$ and on the constant $K$ such that
\begin{equation}
\caG(n\mid n^{\infty})\leq K\,,
\end{equation}
all the details being given in full into \cite{FETA17}. Finally, the term $H_{1}(t)$ in (\ref{C2RD}) is given by:
\begin{equation}
\sum_{j=1}^{s}\bigg[\bigg(\sqrt{\frac{n}{n^{\infty}}}\bigg)^{a^{h}}-\bigg(\sqrt{\frac{n}{n^{\infty}}}\bigg)^{b^{h}}\bigg]^{2}\geq H_{1}(t)\sum_{i=1}^{k}(\sqrt{\frac{n_{i}}{n_{i}^{\infty}}}-1)^2\,.
\end{equation}

\subsection{The Diffusion-Drift approximation}

Let $d\approx m=O(1)$, that is
\begin{equation}
\frac{k_{B}\theta}{e^{2}}\approx\frac{1}{\epsilon L}\,,
\end{equation}
and
\begin{equation}
\frac{k}{d}=\frac{eL^{2}}{k_{B}\theta}\frac{\kappa}{\mu}=o(\varepsilon)\,,\quad\frac{k}{m}=\frac{\epsilon L^{3}}{e}\frac{\kappa}{\mu}=o(\varepsilon)\,,
\end{equation}
then we recover, to within higher order terms, the Poisson-Nernst-Planck system, used to model ion fluxes, cell biology and other electrically driven evolution phenomena:
\begin{eqnarray}\label{BVPNERNST}
&\dive(D[\nabla n]+MN[ez\otimes\nabla\varphi])=\dot{n}\,,&\nonumber\\\
&&\mbox{ in }\Omega\times[0\,,\tau)\,,\\
&-\epsilon\Delta\varphi=ez\cdot n\,,&\nonumber\
\end{eqnarray}
with Neumann boundary conditions. This model was used into \cite{LGW11} and \cite{LGWBW11} to describe the initial stage of the scintillation, before the recombination takes place. The system (\ref{BVPNERNST}) is well-studied and there are many results, dealing with existence, asymptotic decay, equilibrium solutions and even explicit analytical solutions for $k=2$ (\emph{e.g.\/} amongst the many \cite{MAC70}-\cite{DYA05} and for more recent advances \cite{BSW12}-\cite{GAFU21}).

Equation (\ref{BVPNERNST})$_{1}$ can be put in the equivalent gradient flow formulation (\ref{gradflow2}) with $\Kel=\Kel_{D}$. Its stationary solutions $(n^{\infty}\,,\varphi^{\infty})$ are characterized by $s^{\infty}=$const. and therefore (\ref{stationary2}) is replaced by
\begin{equation}\label{stationary3}
n_{j}^{\infty}(x)=c_{j}\exp(\frac{s^{\infty}-ez_{j}\varphi^{\infty}(x)}{\theta k_{B}})\,,\quad j=1\,,2\,,\ldots\,,k\,,
\end{equation}
and the Poisson-Boltzmann equation (\ref{LBG}) changes accordingly.

\subsection{The Reaction-Drift approximation}

I this case $k\approx m=O(1)$, and
\begin{equation}
\frac{d}{k}=\frac{\theta k_{B}}{eL^{2}}\frac{\mu}{\kappa}=o(\varepsilon)\,,\quad\frac{d}{m}=\frac{\epsilon L\theta k_{B}}{e^{2}}=o(\varepsilon)\,.
\end{equation}
The boundary value problem (\ref{BVPadim}) reduces in this case to 
\begin{equation}\label{electrodrift}
\dive(M[N(n)z\otimes\nabla\varphi])-r(n)=\dot{n}\,,\quad\mbox{ in }\Omega\times[0\,,\tau)\,,
\end{equation}
and (\ref{laplace}). At the best of our knowledge this equation was never used for scintillators and indeed it can be found in the theory of dopant diffusion in semiconductors (\emph{vid. e.g.\/} \cite{HFS95}): however we list it together with those used for scintillators for the sake of completeness. For the boundary value problem (\ref{electrodrift}) and (\ref{laplace}) globale existence, uniqueness and asymptotic decay results were obtained into \cite{GLIH97b} (\emph{vid.\/} also \cite{GLIH97} and \cite{GL08}).

\subsection{The Kinetic approximation}

Let $\varepsilon$ a small parameter, then for $k=O(1)$ and
\begin{equation}
\frac{d}{k}=\frac{\theta k_{B}}{eL^{2}}\frac{\mu}{\kappa}=o(\varepsilon)\,,\quad\frac{m}{k}=\frac{e}{\epsilon L^{3}}\frac{\mu}{\kappa}=o(\varepsilon)\,,
\end{equation}
then (\ref{BVPadim})$_{2}$ reduces, to within higher order terms in $\varepsilon$, to the rate ODE equation:
\begin{equation}\label{kinetic1}
-r(n)=\dot{n}\,,\quad\mbox{ in }[0\,,\tau)\,,
\end{equation}
with initial condition $n(0)=n_{o}$; we remark that by (\ref{kinetic1}) and (\ref{firstinte}), then
\begin{equation}
C_{z}=z\cdot n\,,\quad C_{z}\in\Real\,,
\end{equation}
is a first integral of (\ref{kinetic1}). Also in this case (\ref{BVPadim})$_{2}$ is uncoupled, the solution of (\ref{kinetic1}) becoming a charge density supply.

With such approximation, which describes phenomena mainly driven by recombination, we recover the so-called \emph{Kinetic model}, the oldest and most used phenomenological model for scintillators \cite{GBD15}, \cite{BM09}, \cite{BD07}-\cite{VA17},  which is borrowed from the kinetic of chemical reactions, \emph{see e.g.\/} \cite{WA57}, \cite{WA58}.
 
In a series of recent papers, \cite{HMMI16}-\cite{MAMI18}, dealing with the gradient structure of (\ref{RDD}) it is show that also (\ref{kinetic1}) admits a gradient structure with $\Kel=\Kel_{R}$ for detailed balanced recombination mechanisms and hence, the energy dissipation methods can be used with success to get asymptotic decay estimates. 

In detail, into \cite{MAMI18} the well-posedness of (\ref{kinetic1}) was obtained, in the sense that for all initial data $n_{o}\in\caM$ there exists an unique global solution $n:[0\,,\tau)\rightarrow\caM$. Moreover it was proved that since the recombination term satisfies the detailed balance conditions for the equilibrium state $n^{\infty}=c$, then (\ref{kinetic1}) admits the gradient structure:
\begin{equation}
\dot{n}=-H\caD U(n)\,,
\end{equation}
with $H$ given by (\ref{Hmatrix}).

Finally, into \cite{MIMI18} an estimate of the kind of (\ref{estimate0}) was provided:
\begin{equation}
\caD^{*}\geq C_{R}\,\caG(n\mid n^{\infty})\,,
\end{equation}
where the simplified dissipation $\caD^{*}$, such that $\caD\geq k_{o}\caD^{*}$ where $k_{o}$ is given by \emph{e.g.} (\ref{bound0}), is defined as:
\begin{equation}
\caD^{*}=\sum_{h=1}^{s}\frac{k_{h}}{k_{o}}\ell(\frac{n^{a^{h}}}{c^{a^{h}}}\,,\frac{n^{a^{h}}}{c^{b^{h}}})\,.
\end{equation}

\subsection{The Diffusive approximation}

For $d=O(1)$ (which, by (\ref{diffL}) means that $L\approx L_{D}$), and
\begin{equation}
\frac{m}{d}=\frac{e^{2}}{\epsilon L\theta k_{B}}=o(\varepsilon)\,,\quad\frac{k}{d}=\frac{eL^{2}}{\theta k_{B}}\frac{\kappa}{\mu}=o(\varepsilon)\,,
\end{equation}
then from  (\ref{BVPadim})$_{1}$ we obtain a classical anisotropic parabolic equation:
\begin{equation}\label{heatN}
\dive D[\nabla n]=\dot{n}\,.
\end{equation}
In this case, which describes the phenomena when the diffusion within the track is the driving mechanism of scintillation, we recover the \emph{Diffusive model}, which describes the diffusion of excitation carrier within the track  \emph{vid.\/}  \cite{BM12}, and is rarely  used alone as in \emph{e.g.\/} \cite{KD12}. Once again the equation of electrostatic (\ref{BVPadim})$_{2}$ is uncoupled. We notice that the ratio $k/d$ is known as the \emph{Thiele modulus} in the kinetics of chemical reactions. 

A sharp estimate of the decay time for (\ref{heatN}) can be easily found by following \cite{DFM08}:\footnote{We notice that in \cite{HMMI16} this problem was instead formulated in terms of gradient structure with the relative entropy as driving functional, as it was done for the general reaction-diffusion-drift systems.} first of all we notice that at the equilibrium $\dot{n}=0$ and, by the Neumann boundary conditions, $\nabla n=0$, hence $n^{\infty}=\bar{n}^{\infty}=\mathrm{const}$ and by the charge conservation (\ref{conservationcharge}):
\begin{equation}
z\cdot\bar{n}(t)=z\cdot\bar{n}_{o}=z\cdot n^{\infty}=0\,.
\end{equation}
Further, by multiplying both sides of (\ref{heatN}) by $n(x\,,t)-n^{\infty}$, integrating by parts with the Neumann boundary conditions and by the Poincar\'e  inequality then we get:
\begin{equation}\label{poinac1}
\frac{\dif}{\dif t}\int_{\Omega}|n(x\,,t)-n^{\infty}|^{2}=-2\int_{\Omega}D[\nabla n]\cdot\nabla n\leq -\frac{2\delta}{\caL(\Omega)}\int_{\Omega}|n(x\,,t)-n^{\infty}|^{2}\,.
\end{equation}
The first and last terms of (\ref{poinac1}) yield a first order equation which can be integrated to obtain at an expression alike (\ref{estimate0}):
\begin{equation}
\|n(x\,,t)-n^{\infty}\|_{L^{1}(\Omega)}^{2}\leq\caG(n_{o}\,,n^{\infty})\exp(-\frac{2\delta}{\caL(\Omega)}t)\,,
\end{equation} 
with
\begin{equation}
\caG(n_{o}\,,n^{\infty})=\int_{\Omega}|n_{o}(x)-n^{\infty}|^{2}\,.
\end{equation}

Accordingly the decay time $\tau_{d}$ can be estimated as:
\begin{equation}
\tau_{d}\leq\frac{\caL(\Omega)}{2\delta}=\frac{e\caL(\Omega)}{2k_{B}\theta\mu}\,.
\end{equation}

Whenever $D=\diag\{D_{1}\,,D_{2}\,,\ldots\,,D_{k}\}$, equation (\ref{heatN}) reduces to the $k$ independent classical diffusion equations:
 \begin{equation}\label{heat}
D_{j}\Delta n_{j}=\dot{n}_{j}\,,\quad\mbox{ in }\Omega\times[0\,,\tau)\,, \quad D_{j}>0\,,j=1\,,2\,,\ldots\,,k\,,
\end{equation}
with Neumann boundary conditions, whose a complete mathematical treatment can be found into many books, like \emph{e.g.\/} \cite{EDB}.

\section{Conclusions}

For the reaction-diffusion-drift equation which describes the evolution and recombination processes of charge carriers in scintillators we gave an overview of the existence and asymptotic decay estimate which are know to date, at the best of our knowledge. Despite the fact that the topics is a well-studied one, as the non-exhaustive list of references shows, there are still many unanswered questions which deserves further investigations: here we shall give a concise list of some them. Such a list is of course far from be exhaustive and its items are those which seem more interesting for the physicists.

First of all the mathematical treatment of the RDD system is based on a precise choice of the entropic term $F(n)$, namely the Gibbs entropy: however the physics of scintillation processes suggests that such a choice can be appropriate into describing the behavior over a limited range of energy. The choice of Fermi-Dirac potential in place of the Gibbs entropy will be not only more general but also more related to the true physical nature of the phenomena. Clearly, this leads to a different formulation for the recombination term which requires a different mathematical treatment.

A second point of interest is the decay time estimate obtained for the RDD system: here there are two major directions for further investigations. The first one is to extend the results obtained into \cite{FK18} for a system of two charge carriers, $k=2$ to systems with general $k$ as most of the scintillator phenomenological studies require. The second and more intriguing aspect is the following: the decay time estimate which follows from the result of \cite{FK18} gives as an upper bound of the  decay time the maximum between two values, one which depends on the carrier mobility and the other which depends on the recombination time. When into \cite{D3} we applied these results to four scintillating crystal, not only we obtained a very good estimate ad consistent  of the \emph{fast} decay time: we also notice that the \emph{minumum} between these two values looked like a \emph{lower bound} for the {slow} decay time. Clearly the results of \cite{FK18} tells nothing about this but it would be tempting to prove such and assertion. In such a case we should have two bounds, an upper one on the fast decay time and a lower one on the slow decay time, a result that would be appreciable in terms of material science. 

A third remark is that the \emph{Entropic methods} used to get these estimate should be applied also to the approximated phenomenological models, in particular to the \emph{kinetic} one: we need to remark that a very useful aspects of the results in \cite{FK18} is the explicit dependence of the estimates on the constitutive parameters of the RDD systems. This should allows for a predictive use of these results, not only for an \emph{a-posteriori}  check with the available experimental data.

As far as the questions related to the existence issues for very general form of recombination term is concerned, the necessity of further generalization doesn't need to be justified, since they are the foundation stone of any numerical procedure we need to implement in order to obtain numerical solutions.

Finally there are two side-aspects related to those we looked at in this paper and which we left out. The first one concerns the solution of the Bethe-Bloch equation, which is a necessary requirement to bridge the microscopic world with the phenomenological treatment we want to give to scintillation. To date there are many results which concern mainly the radiation decay length or the Bragg peak: a treatment which conveys in a straightforward way the relevant parameters of the phenomena at the mesoscopic scale would be welcomed.

The second and in some sense more important problem we left out is the \emph{Ligth yield}: such a parameter is indeed a measure of scintillator effectiveness and tells us about the minimum energy we can detect. To date there is not a coherent mathematical definition for this parameter, apart a descriptive one:
\begin{equation*}
LY=\frac{\mbox{Number of charge carriers recombined into photons}}{\mbox{Total number of charge carriers generated}}\,,
\end{equation*}
which in turns can be applied either locally or globally, whatever this means. Once again, a formal definition and its consequent mathematical treatment are, in our opinion, still missing and it would benefit from the huge amount of experimental works which give a precise evaluation of Light Yield.

\section*{Acknowledgements}

\noindent The research leading to these results is within the scope of CERN R\&D Experiment 18 "Crystal Clear Collaboration" and the PANDA Collaboration at GSI-Darmstadt. I wish to thanks Nella Rotundo for her many suggestions, Klemens Fellner, Michael Kniely and Bao Quoc Tang for giving me the possibility to talk about the mathematical aspects of scintillation in a lively, despite remote, seminar and finally Paolo Maria Mariano for some useful comments.

\end{document}